# ЖАРАТЫЛЫСТАНУ ЖӘНЕ ТЕХНИКАЛЫҚ ҒЫЛЫМДАР

# NATURAL AND TECHNICAL SCIENCES

*Batyr Sharimbayev[1*], Shirali Kadyrov[2], Aleksei Kavokin[3]*
[1,3]SDU University, Kaskelen, Kazakhstan
[2]New Uzbekistan University, Tashkent, Uzbekistan
*e-mail: batyr.sharimbayev@sdu.edu.kz

## DEVELOPMENT AND OPTIMIZATION OF PHYSICS-INFORMED NEURAL NETWORKS FOR SOLVING PARTIAL DIFFERENTIAL EQUATIONS

**Abstract:** This work compares the advantages and limitations of the Finite Difference Method with Physics-Informed Neural Networks, showing where each can best be applied for different problem scenarios. Analysis on the $L_2$ relative error based on one-dimensional and two-dimensional Poisson equations suggests that FDM gives far more accurate results with a relative error of $7.26 \times 10^{-8}$ and $2.21 \times 10^{-4}$, respectively, in comparison with PINNs, with an error of $5.63 \times 10^{-6}$ and $6.01 \times 10^{-3}$ accordingly. Besides forward problems, PINN is realized also for forward-inverse problems which reflect its ability to predict source term after its sufficient training. Visualization of the solution underlines different methodologies adopted by FDM and PINNs, yielding useful insights into their performance and applicability.

***

**Аннотация**: В этой работе сравниваются преимущества и ограничения метода конечных разностей с использованием нейронных сетей, основанных на физике, и показывается, где каждый из них может быть наилучшим образом применен для различных сценариев решения задач. Анализ относительной погрешности $L_2$, основанный на одномерных и двумерных уравнениях Пуассона, показывает, что FDM дает гораздо более точные результаты с относительной погрешностью $7,26 \times 10^{-8}$ и $2,21 \times 10^{-4}$, соответственно, по сравнению с PINNs, с погрешностью $5,63 \times 10^{-6}$ и $6,01 \times 10^{-3}$ соответственно. Помимо прямых задач, PINN реализован также для прямых и обратных задач, что отражает его способность прогнозировать исходный код после достаточного обучения. Визуализация решения подчеркивает различия в методологиях, принятых FDM и PINNs, и дает полезную информацию об их производительности и применимости.

**Ключевые слова:** численный анализ, прямые и обратные задачи, глубокое обучение, PINN, FDM.

***

**Аңдатпа:** Бұл жұмыс Ақырлы Айырмашылық Әдісінің артықшылықтары мен шектеулерін Физикаға Негізделген Нейрондық Желілермен салыстырады, олардың әрқайсысын әртүрлі проблемалық сценарийлер үшін қай жерде қолдануға болатынын көрсетеді. Бір өлшемді Және екі өлшемді Пуассон теңдеулеріне негізделген $L_2$ салыстырмалы қатесін талдау FDM Салыстырмалы қателігі сәйкесінше $7{,}26 \times 10^{-8}$ және $2{,}21 \times 10^{-4}$ Болатын Әлдеқайда дәл нәтижелер беретінін көрсетеді.Түйреуіштермен салыстырғанда, қателігі сәйкесінше $5{,}63 \times 10^{-6}$ және $6{,}01 \times 10^{-3}$. Сонымен қатар, алға бағытталған проблемалар, PINN жеткілікті дайындықтан кейін бастапқы мерзімді болжау қабілетін көрсететін алға-артқа есептер үшін де жүзеге асырылады. FDM және PINNs қабылдаған әртүрлі әдістемелердің шешімдерін суреттеу, олардың өнімділігі мен қолданылуы туралы пайдалы түсінік береді.

**Түйін сөздер:** сандық талдау, тура-кері есептер, терең оқыту, PINN, FDM.

*1. Introduction*

PDEs would form a crucial backbone in the understanding and modeling of various real-world problems. In simple words, PDEs describe how something changes through time and space based on very well-defined mathematical rules. The central issue with PDEs is determining whether these rules suffice to guarantee a unique solution to the problem [1][2].

In most real-life cases, it is impossible to find the exact solutions for such complex equations involving PDEs. This gave rise to the development of numerical methods for estimating solutions. Some of the popular methods include the finite element method, finite difference method, finite volume method, and spectral element method. Out of these, FEM is the most advanced with strong mathematical support for ensuring accurate results, stability, and error control. FEM solvers often employ efficient techniques such as sparse linear systems or iterative methods, which make them suitable for many practical problems [3].

Most recently, new techniques have emerged due to the growth of deep learning which aid in the resolution of PDEs. One of these is Physics-Informed Neural Networks (PINNs), where the solution to the PDE is parameterized as a neural network that is trained to minimize the residual errors of the PDE along with its boundary conditions. The use of PINNs has been beneficial for problems like grid dependency and dealing with high dimensional problems. They also happen to be tremendously useful for solving various kinds of PDEs because of their simplicity in implementation and direct incorporation of the underlying physics within the learning itself [4].

Despite their adoption, deep learning algorithms, including PINNs, still suffer from high computational costs, complex optimization, and weak theoretical foundations. The purpose of this study is to carry out a computer-based assessment in an attempt to test the usability and effectiveness of the finite element method against the recent Physics Informed Neural Networks approach. To do this, we first need to discuss relevant developments and implementations of PINNs in the context of PDE solving.

1.1 Literature review

PINNs have recently gained significant attention for solving complex problems related to PDEs. In one work [5], the authors presented an approach involving the use of PINNs for modeling wave propagation and performing full waveform inversions. This was applied on the acoustic wave equation and showed good performance on a suite of different synthetic scenarios that included complex structures. Unlike the classical approaches, such as finite difference or spectral element methods, which rely on grids, PINNs are "meshless". They can handle arbitrary boundary conditions, including absorbing boundaries, without additional modifications. While classical methods are still much faster and more accurate for forward simulations, PINNs have demonstrated quite impressive capabilities in solving inversion problems, including rather complicated ones. The present study has introduced PINNs as a flexible and promising tool in solving seismic problems, as well as any other application that includes the solution of PDEs.

On the basis of PINNs, another work [6] introduced a different method called Variational Physics-Informed Neural Networks. VPINNs extended the original framework of PINNs by incorporating the variational form of the PDE in its training process. Instead of relying on randomly selected points in the domain, VPINNs employed a mixture of neural networks ("trial space") and mathematical polynomials ("test space") to achieve higher accuracy. This approach reduces the cost of training and gives more precise results, especially with simpler neural networks. It was shown that VPINNs performed much faster and more accurately than PINNs; hence, they may be a suitable alternative in solving nonlinear PDEs.

Another work [7] compared PINNs with the finite element method -a classic solution technique for PDEs. Several types of PDEs were considered: Poisson's equation, Allen-Cahn equation, and semilinear Schrödinger equation in different dimensions. The focus was on computational cost and the accuracy of the two methods. The current paper determines that each approach has its strong points; FEM may turn out to be more reliable in some problems, while in others, PINNs allow great flexibility and ease of use.

2. Mathematical Background
2.1 Poisson Equation

The Poisson equation is one of the most basic PDEs, which results in a potential field created by a given source. It provides significant applications to physics and engineering, modeling electrostatics, heat conduction, and fluid dynamics among many others. References [8][9][10].

In one dimension, the Poisson equation is given by:

$$\frac{d^2 u(x)}{dx^2} = f(x)$$

where $u(x)$ is the unknown function, $f(x)$ is a source term.

In such cases, the domain is normally an interval $[a, b]$ and the boundary conditions are often given at both ends, say in the following form of Dirichlet boundary conditions::

$$u(a) = u_0, \quad u(b) = u_1$$

In two dimensions, the Poisson equation generalizes to:

$$\nabla^2 u(x, y) = f(x, y)$$

where is the 2D Laplacian operator, represents the unknown function, and is called the source term. In two dimensions, it can normally be a rectangular region or even more complicated geometry. On edges, boundary conditions are imposed. There may also be Dirichlet or Neumann boundary conditions on these. A typical Dirichlet boundary condition takes the following form: $u(x, y) = g(x, y)$ on the boundary of the domain.

The accuracy of the solution approximations for the 1D and 2D Poisson equations is evaluated using the $L_2$ norm and relative error [14][15]. The $L_2$ norm of a vector $v = [v_1, v_2, \ldots, v_n]$ is defined as:

$$\|v\|_2 = \sqrt{\sum_{i=1}^{n} v_i^2}.$$

This norm gives the size of the magnitude of a vector in Euclidean space. To verify these results for accuracy, we calculate the $L_2$ relative error between the computed solution $\hat{u}$ and the exact solution u defined as:

$$L_2 = \frac{\|\hat{u} - u\|_2}{\|u\|_2}.$$

*2.2 Physics-Informed Neural Networks (PINN)*

Let $u_\theta(x, y)$ denote the neural network approximation of $u(x, y)$ where θ is used to denote the parameters of the network. The network architecture considered here is a fully connected feed-forward neural network, such that the application of activation functions is performed layer by layer. It can be written as:

$$u_\theta(x, y) = \mathrm{NN}_\theta(x, y)$$

The network layers are defined as follows:

$$x^{(i+1)} = \sigma\left(w^{(i)} x^{(i)} + b^{(i)}\right)$$

for each hidden layer, and the output layer is linear. The residual of the PDE is defined as:

$$R(x, y) = \frac{\partial^2 u_\theta}{\partial x^2} + \frac{\partial^2 u_\theta}{\partial y^2} - f(x, y)$$

where $f(x, y)$ is the source term. The PDE loss is computed as the mean squared error of the residual over a set of domain points $(x_i, y_i)$:

$$L_{\mathrm{PDE}} = \frac{1}{N} \sum_{i=1}^{N} R(x_i, y_i)^2$$

The boundary conditions loss is computed similarly, where the loss for boundary points $(x_i, y_i)$ is given by:

$$L_{\mathrm{BC}} = \frac{1}{M} \sum_{j=1}^{M} (u_\theta(x_j, y_j) - u(x_j, y_j))^2$$

The total loss function combines the PDE residual loss and the boundary conditions loss:

$$L(\theta) = L_{\mathrm{PDE}} + L_{\mathrm{BC}}$$

### 2.3 Finite Difference Method (FDM)

The FDM is a numerical method for the solution of PDEs that is based on discretizing their solutions on a lattice of points that discretize the domain. Based on this concept, the FDM relies upon approximating the derivatives of an unknown function in terms of finite differences that replace the PDE by a set of algebraic equations [11][12][13].

First, consider the one dimensional Poisson equation:

$$\frac{d^2 u(x)}{dx^2} = f(x)$$

Using a finite difference scheme, the second derivative of can be approximated by the following way :

$$\frac{d^2 u(x)}{dx^2} \approx \frac{u(x+\Delta x) - 2u(x) + u(x-\Delta x)}{(\Delta x)^2}$$

where $\Delta x$ is the step size and $u(x)$ is the unknown function at the grid points. By discretizing the domain $[a, b]$ with a grid of points, we convert the continuous PDE into a system of algebraic equations which can be solved numerically.

In two dimensions, the Poisson equation is given by:

$$\nabla^2 u(x, y) = f(x, y)$$

where $\nabla^2$ is the Laplacian operator:

$$\nabla^2 = \frac{\partial^2}{\partial x^2} + \frac{\partial^2}{\partial y^2}.$$

The second derivatives in the $x$ and $y$ directions are approximated by finite differences. The Laplacian operator in 2D is approximated as:

$$\nabla^2 u(x, y) \approx \frac{u(x+\Delta x, y) + u(x-\Delta x, y) + u(x, y+\Delta y) + u(x, y-\Delta y) - 4u(x, y)}{(\Delta x)^2}$$

where $\Delta x$ and $\Delta y$ are the step sizes in the $x$ and $y$ directions, respectively.

The FDM discretizes the domain of both 1D and 2D Poisson equations, then solves the resulting system of algebraic equations. The solution is carried out by an iteration over all the grid points, with the values of $u(x)$ or $u(x, y)$ updated in conformity with the finite difference approximation delivered by the PDE. The boundary conditions are normally prescribed along the edges of the domain in terms of either dirichlet or neumann-type conditions.

### 3. Methodology

This work concentrates on developing and enhancing PINNs for the solutions of the 1D and 2D Poisson equations and comparing them with the FDM.

The data to be used in these tasks were sampled from a given range by using the Latin Hypercube Sampling (LHS) technique. The range of interest for the problem at hand was between 0 and 1, and the LHS technique was used in order to draw 256 samples from the range and 2 from the boundary conditions. The obtained samples were used in training the PINN as well as testing the prediction capabilities of the model.

To begin with, FDM is the first approach we had to study to solve the Poisson equation. The domain of the 1D example is split into a lattice grid, and the equation is solved through second derivation approximation. The same procedure is applied in 2D, in which both spatial dimensions are discretized. The numerical solution is obtained and compared to the calculated exact solutions in order to check accuracy using the $L_2$ relative error.

At this stage, the PINN is implemented to solve the equations. The neural network is trained by minimizing a loss function which comprises the Poisson equation and the bounding conditions. The training is therefore carried out in two steps: firstly, with the Adam optimizer to fit the parameters of the model, and subsequently with the L-BFGS method for fine tuning.

Latin Hypercube Sampling, in both 1D and 2D scenarios, provides points for model analyses. The network is trained to provide a set of boundary conditions. The performance of the model is determined by comparing the output of the neural network with a known solution using the relative $L_2$ norm.

Ultimately, we applied the FDM to a thermal problem involving a second order ordinary differential equation (ODE). We split the spatial domain into parts, provided intervals, and boundary conditions to both ends. The problem was solved iteratively until the difference between subsequent values was less than some small predefined threshold. The materials' properties and source term were represented as functions of spatial position. The PINN solution was compared to the one obtained using multivariate interpolation to assess its accuracy. Then the forward-inverse problem was solved using a PINN structure. The PINN was based on a feedforward neural network with hidden layers that were trained to estimate temperature and source term via a learned loss function. The model was validated against FDM and observational data. The code was written with the python and is available on github: https://github.com/hardkazakh/pinn-vs-fdm

### 4. Experiential Results
### 4.1 1D Poisson equation

Let us investigate 1D Poisson equation defined as:

$$\frac{d^2 u(x)}{dx^2} = 16x^7 e^{-x^4} - 20x^3 e^{-x^4}, \quad x \in [0, 1], \quad (3.1)$$

with Dirichlet boundary conditions:

$$u(0) = 0, \quad u(1) = e^{-1}.$$

The equation has an analytical solution, which can be written as:

$$u(x) = xe^{-x^4}.$$

In Section 2.3, we explained that the first step to solving PDEs with the

FDM is to rewrite the equation in a weak form. We already did this for the Poisson equation. The next step is to create a mesh. This is like breaking the interval [0, 1] into small pieces, called cells. The number of cells is 512. More cells mean the grid is finer, which gives a more accurate solution, but it also takes more time and computing power.

For solving the 1D Poisson equation using PINNs, there are three design parameters that we need to specify before training. The first step is choosing a loss function. Following the vanilla PINNs approach, we evaluate the goodness of the solution using the discretised mean squared error over the PDE, boundary and initial conditions. In particular, we define the loss function as:

$$\text{Loss}(\theta) = \frac{1}{N} \sum_{i=1}^{N} \|\Delta u_\theta(x_i) - (16x^7 \exp(-x^4) + 20x^3 \exp(-x^4))\|_2^2 + \|u_\theta(0)\|_2^2 + \|u_\theta(1) - \exp(-1)\|_2^2$$

with $u_\theta$ the neural network, $\theta$ the trained weights, and 512 of collocation points $x_i$ sampled in each epoch using LHS.

The second design parameter is the neural network architecture, that is, the type of neural network, the activation function, and the number of hidden layers and nodes. For the 1D Poisson case, we train feed-forward dense neural networks with *tanh* as the activation function. We use the result on the architecture of [20, 20, 20, 1].

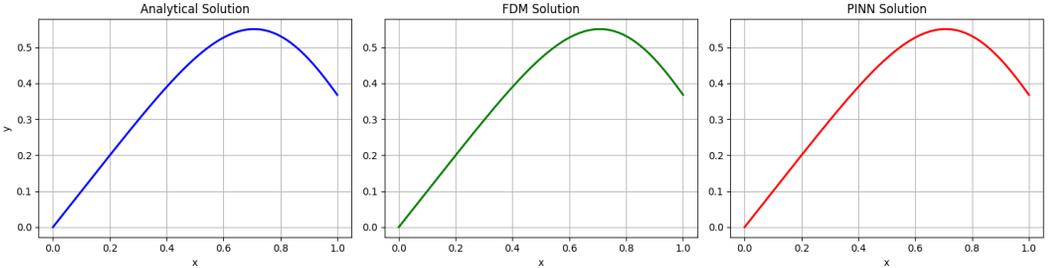

Figure 1. Comparison of solutions: 1-Exact, 2-FDM, and 3-PINN.

The approximations of the 1D Poisson equation solution using FEM and PINNs are compared to the exact solution on a [0, 1] interval with 512 points. Figure 1 shows the exact solution and the approximations. One PINN setup, with only one hidden layer and one node, performs poorly and fails to satisfy the boundary conditions. All approximations are very close to the exact solution.

For the 1D Poisson equation, the relative error for the FDM is calculated to be $7.26 \times 10^{-8}$, while the relative error for the PINN approach is $5.63 \times 10^{-6}$. These results show that FDM provides a more accurate approximation of the solution compared to PINNs in the 1D case.

*4.2 2D Poisson equation*
1) Let us now investigate 2D Poisson equation defined as:

$$\frac{\partial^2 u(x,y)}{\partial x^2} + \frac{\partial^2 u(x,y)}{\partial y^2} = 2\left(x^4(3y-2) + x^3(4-6y) + x^2\right). \quad (3.2)$$

The boundary conditions are:

$$u(x,0) = u(x,1) = u(0,y) = u(1,y) = 0.$$

The analytical solution of the equation is:

$$u(x,y) = (x-1)^2 y(y-1)^2 x^2.$$

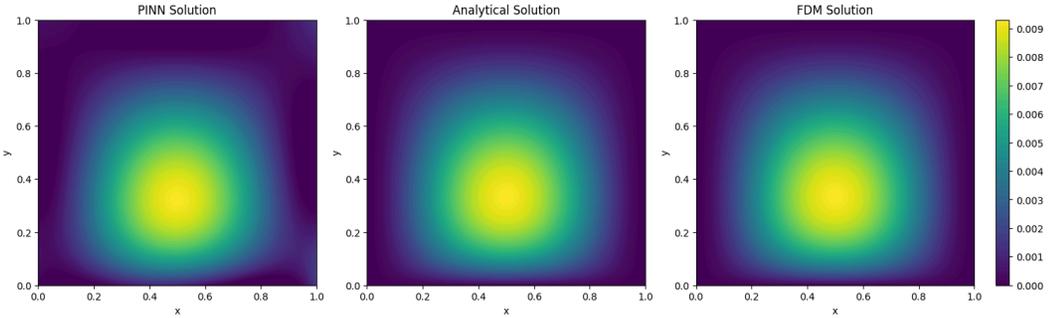

Figure 2. Comparison of solutions: 1-PINN, 2-Exact, and 3-FDM.

Figure 2 shows the analytical and approximate solutions of the 2D Poisson equation. For the 2D Poisson equation, the $L_2$ relative error for the FDM approximation is $2.21 \times 10^{-4}$, while for the PINN, it is $6.01 \times 10^{-3}$. Again, the FDM method shows a significantly lower error compared to the PINN, indicating that FDM achieves a more precise solution in both 1D and 2D cases.

These error analyses again confirm the accuracy of the FDM approach in solving Poisson equations, especially in comparison with PINNs, which showed higher relative errors in both 1D and 2D problems. However, in higher order equations, PINN can give better results. We consider this in our future research.

*4.3 PINN approach for simultaneous forward-inverse problems (FIP)*

PINN embeds both forward and inverse problems under one framework by embedding the data and physical laws in the loss function of the neural network. The network jointly predicts the forward solution of a PDE and estimates unknown parameters or inputs (inverse problem) through the optimization of one combined loss function that considers residual of PDE, boundary conditions, and discrepancies between model predictions and observations. This capability of handling both tasks together makes PINNs highly effective at solving problems that involve solution estimation and identification of parameters.

The problem is defined by the following second-order differential equation:

$$\frac{d^2 U(x)}{dx^2} - a(x)U(x) = Q(x), \quad 0 < x < L, \tag{3.3}$$

where $Q(x)$ represents the source term and $a(x)$ describes the varying coefficient. Specifically, the source term is given by:

$$Q(x) = 1 + b_1 \sin(w_1 x)$$

and the varying coefficient is defined as:

$$a(x) = b_1 + \frac{x}{1+x^2}.$$

The boundary conditions for the differential equation are given by: $U(0) = 1$, $U(L) = 3$. The problem now involves the solution of this second-order differential equation along with the boundary conditions shown above. $Q(x)$, the source function, incorporates a sine function that could model some periodic influence in the system. This coefficient, $a(x)$, depends on the position x in the domain; therefore, this makes the equation more complicated, introducing spatial dependence in the solution. Figure 3: Visualization of FDM approximation and source term.

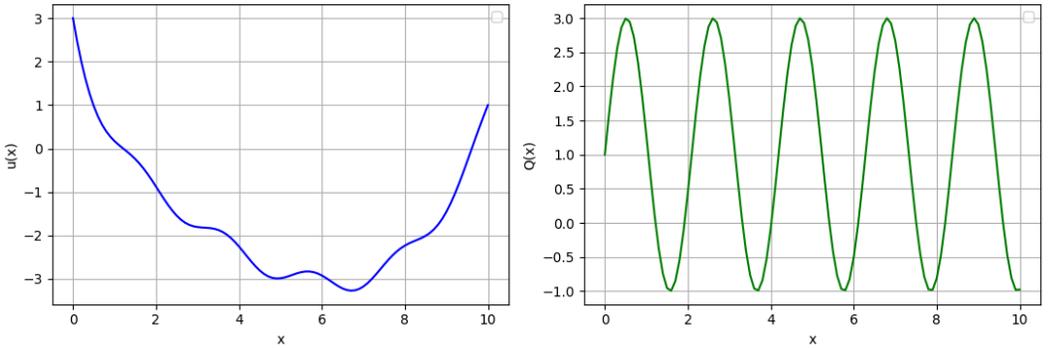

Figure 3. Visualization of FDM Solution and source term $Q(x)$.

The aim of the Experiment 1 is to approximate both the source term $Q(x)$ and the solution $U(x)$ of the given differential equation. PINN simultaneously predicts the solution $U(x)$ while reconstructing the source term $Q(x)$ using the varying coefficient $a(x)$ as part of the system.

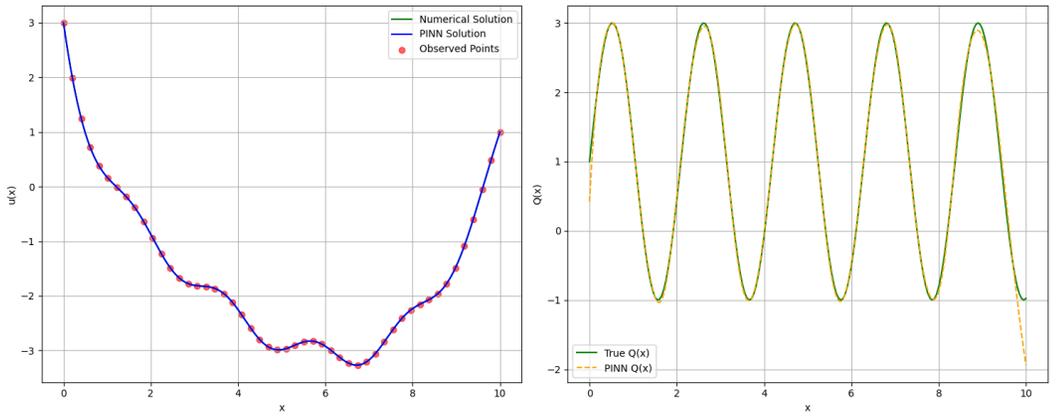

Figure 4. Visualisation of PINN Solution and prediction of $Q(x)$.

The PINN model was trained for 40,000 epochs. Initially, the loss function was quite large, but it gradually decreased as the model learned. By the end of the training, the loss had reduced to $1.87 \times 10^{-2}$, indicating that the model had learned the underlying physics of the problem. The predictions made by the PINN for the temperature distribution and source term are closely aligned with the FDM results. These predictions were accurate not only at the observed points but also at unobserved locations, highlighting the PINN's ability to generalize well across the entire domain. Figure 4 presents the visualization of the PINN's predicted temperature distribution and the corresponding predicted source term.

The aim of Experiment 2 is to approximate both the varying coefficient $a(x)$ and the solution $U(x)$. By incorporating $a(x)$ as an unknown parameter in the system, the model predicts the solution $U(x)$ while reconstructing $a(x)$ from the given data.

The PINN model was trained for 40,000 epochs. Initially, the loss function exhibited high values, but with training, it gradually decreased as the model captured the intricate relationships within the system. By the end of the training, the loss had reduced to $3.1563 \times 10^{-2}$, demonstrating that the model effectively learned both the solution and the coefficient .

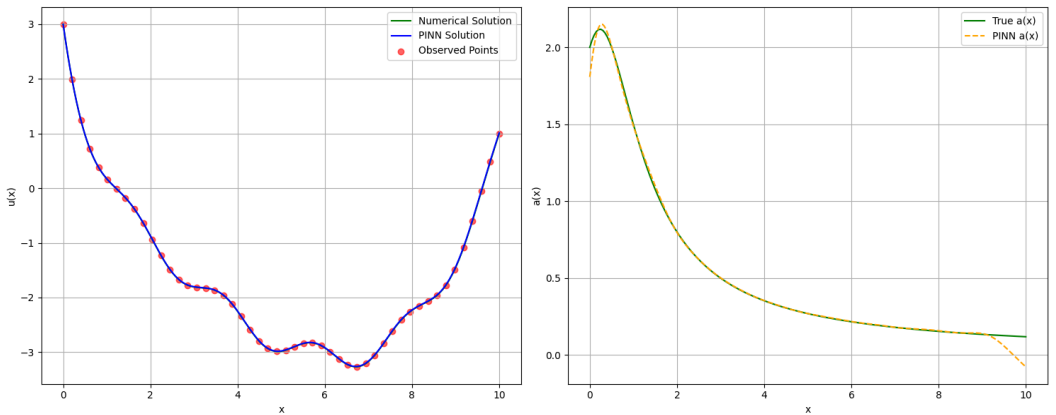

Figure 5. Visualisation of PINN Solution and prediction of $a(x)$.

Figure 5 illustrates the visualization of the PINN's predictions for the temperature distribution and the reconstructed varying coefficient. The predicted solution from the PINN closely aligns with the numerical solution obtained using FDM. The observed points are well-matched with the PINN-predicted solution, indicating excellent accuracy at sampled locations. Additionally, the smoothness of the PINN solution demonstrates its capability to generalize effectively across unobserved regions.

The PINN accurately reconstructs the varying coefficient, as shown by the close agreement between the true and the PINN-predicted. While small discrepancies may appear at some points, the overall reconstruction is highly accurate, reflecting the model's ability to infer the hidden system parameter.

## 5. Conclusion & Future work

This work compared the accuracy of FDM and PINNs for solving 1D and 2D Poisson equations. Results show that the FDM has much smaller relative errors compared with PINNs for both cases: the relative errors are $7.26 \times 10^{-8}$ versus $5.63 \times 10^{-6}$ for the 1D problem and $2.21 \times 10^{-4}$ versus $6.01 \times 10^{-3}$ for the 2D problem. Although PINNs are less precise for these particular cases, their potential to solve the simultaneous forward-inverse problems showed satisfactory results. The FDM method presented convergence within the allowed iterations, where the obtained results for the temperature distribution and the source term are accurate and in good agreement with the interpolated analytical solution, confirming the reliability of FDM.

For the PINN approach, FIP results for the source term and varying coefficient solutions were:
- FIP Source term: Epoch 40,000, Loss = $1.87 \times 10^{-2}$.
- FIP Varying coefficient: Epoch 40,000, Loss = $3.1563 \times 10^{-2}$.

The temperature distribution, varying coefficient, and source term predictions of the PINN were in excellent agreement with the FDM results, indicating that the model was able to generalize well over the domain, including points that were not observed.

In summary, the results show that both methods work very well: FDM has been proved to produce high-accuracy outcomes, and PINNs have presented a promising alternative in solving both forward and inverse problems in computational physics. Future work will try to further improve the accuracy of PINNs and explore hybrid methods combining the strengths of FDM and PINNs for more complex problems.

## Acknowledgments

This research was funded by the Ministry of Science and Higher Education of the Republic of Kazakhstan within the framework of project AP23487777.